\documentclass[12pt,a4paper]{article}
\usepackage{amssymb,amsmath}
\usepackage{amsfonts}
\usepackage{a4}
\begin{document}
 \baselineskip18pt
 \title{\bf $(1-2u^3)$-constacyclic codes and quadratic residue codes over  $\mathbb{F}_{p}[u]/\langle u^4-u\rangle$}
\author{Madhu Raka, Leetika Kathuria and Mokshi Goyal
 \\ \small{\em Centre for Advanced Study in Mathematics}\\
\small{\em Panjab University, Chandigarh-160014, INDIA}\\
\date{}}
\maketitle
 {\abstract{Let $\mathcal{R}=\mathbb{F}_{p}+u\mathbb{F}_{p}+u^2\mathbb{F}_{p}+u^3\mathbb{F}_{p}$ with $u^4=u$ be a finite non-chain ring, where $p$ is a prime congruent to $1$ modulo $3$. In this paper we study $(1-2u^3)$-constacyclic codes over the ring $\mathcal{R}$, their equivalence to cyclic codes  and find their Gray images. To illustrate this, examples of $(1-2u^3)$-constacyclic codes of lengths $2^m$ for $p=7$ and of  lengths $3^m$ for $p=19$ are given. We also discuss quadratic residue codes over the ring $\mathcal{R}$ and their extensions. A Gray map from $\mathcal{R}$ to $\mathbb{F}_{p}^4$  is defined which preserves  self duality and gives self-dual and formally self-dual codes over $\mathbb{F}_{p}$ from extended quadratic residue codes. }\vspace{2mm}\\{\bf MSC} : 11T71, 94B15.\\
 {\bf \it Keywords }: Constacyclic codes; Negacyclic codes; Self-dual and self-orthogonal codes; Quadratic residue codes, Extended QR-codes, Gray map.}
 \section{ Introduction.}
  $~~~$The class of constacyclic codes plays a significant role in the theory of error correcting codes. These include cyclic and negacyclic codes, which have been well studied since 1950's. Constacyclic codes can be efficiently encoded using simple shift registers. They have rich algebraic structures for efficient error detection and correction, which explains their preferred role in engineering.

Constacyclic codes over finite fields have been investigated by various authors [2,3,6,7,9,10,17] whereas  codes over finite rings have recently generated a lot of interest. Many authors such as [8] worked on constacyclic codes over finite chain rings. In [19], Zhu and Wang investigated $(1-2u)$-constacyclic codes over the finite  ring
$\mathbb{F}_{p}+u\mathbb{F}_{p}$ where $u^2=u$ and defined a Gray map on it. Mostafanasab and Karimi
[15] investigated $(1-2u^2)$-constacyclic codes over the finite ring
$\mathbb{F}_{p}+u\mathbb{F}_{p}+u^2\mathbb{F}_{p}$ where $u^3=u$ and $p$ is an odd prime. Bayram and Siap [4,5] studied constacyclic codes over the ring $\mathbb{F}_{p}[u]/\langle u^p-u\rangle$, where $p$ is a prime. In this paper we extend  results of [15] and discuss $(1-2u^3)$-constacyclic codes over the finite ring
$\mathcal{R}=\mathbb{F}_{p}+u\mathbb{F}_{p}+u^2\mathbb{F}_{p}+u^3\mathbb{F}_{p}$ where $u^4=u$, $p$ is a prime and $p\equiv 1 ({\rm mod}~ 3)$. The unit $(1-2u^3)$ is a special unit of $\mathcal{R}$ as it is inverse of itself. We explore conditions on $n$ and $p$ under which a $(1-2u^3)$-constacyclic code is equivalent to a cyclic code.\vspace{2mm}

Quadratic residue codes are a special kind of cyclic codes of prime length introduced to construct self-dual codes. Quadratic residue codes over $\mathbb{Z}_4$ were studied by Pless and Qian [16].  Kaya  et al. [11] and Zhang  et al [19] studied quadratic residue codes over $\mathbb{F}_p+u\mathbb{F}_p$. Kaya et al [12] studied quadratic residue codes over  $\mathbb{F}_2+u\mathbb{F}_2+u^2\mathbb{F}_2$ whereas Liu et al [13] studied them over non-local ring $\mathbb{F}_p+u\mathbb{F}_p+u^2\mathbb{F}_p$ where $p$ is an odd prime. In this paper we extend their results and discuss quadratic residue codes over the ring
$\mathcal{R}=\mathbb{F}_{p}+u\mathbb{F}_{p}+u^2\mathbb{F}_{p}+u^3\mathbb{F}_{p}$ where $u^4=u$ and $p\equiv 1 ({\rm mod ~} 3)$. We also define a Gray map preserving the self duality  which gives self-dual and formally self-dual codes over $\mathbb{F}_{p}$. In a subsequent  paper we will discuss quadratic residue codes over
$\mathbb{F}_{p}+u\mathbb{F}_{p}+\cdots +u^{m-1}\mathbb{F}_{p}$ where $u^m=u$ and $p\equiv 1 ({\rm mod ~} (m-1))$ for any natural number $m\geq 2$.\vspace{2mm}

The paper is organized as follows: In section 2, we give some basic preliminary results. In Section 3, we  study $(1-2u^3)$-constacyclic codes over the ring $\mathcal{R}$, their equivalence and find their Gray images. To illustrate this, examples of $(1-2u^3)$-constacyclic codes of lengths $5$ and $2^m$ for $p=7$ and examples of $(1-2u^3)$-constacyclic codes of lengths $3^m$ for $p=19$ are given. In Section 4, we study quadratic residue codes over $\mathcal{R}$ and give some of their properties. In Section 5, we study extended quadratic residues codes leading to self-dual and formally self-dual codes over $\mathbb{F}_p$ under another Gray map.
\section{Preliminary results}
 Throughout the paper, $\mathcal{R}$ denotes the commutative ring $\mathbb{F}_{p}+u\mathbb{F}_{p}+u^2\mathbb{F}_{p}+u^3\mathbb{F}_{p}$ where $u^4=u$,
  $p$ is a prime and $p\equiv 1({\rm mod}~3)$. $\mathcal{R}$ is a ring of size $p^4$ and characteristic $p$. For   a primitive element $\alpha$ of $\mathbb{F}_{p}$, take $\xi=\alpha^{\frac{p-1}{3}}$, so that $\xi^3=1, \xi\neq 1$ and $\xi^2+\xi+1=0$. Let $\eta_1, \eta_2, \eta_3,
 \eta_4$ denote the following elements of $\mathcal{R}$. \vspace{2mm}
\begin{equation} \begin{array}{ll}\eta_1=1-u^3\\
\eta_2= 3^{-1}(u+u^2+u^3)\\
\eta_3= 3^{-1}(\xi u+\xi^2u^2+u^3)\\
\eta_4= 3^{-1}(\xi^2u+\xi u^2+u^3)\end{array}\end{equation}\vspace{2mm}

A simple calculation shows that $\eta_i^2=\eta_i$, $\eta_i\eta_j=0$ for $1\leq i, j \leq 4, ~i\neq j$ and $\sum_{i=1}^4 \eta_i=1$. The decomposition theorem of ring theory tells us that  $\mathcal{R}=\eta_1\mathcal{R}\oplus\eta_2\mathcal{R}\oplus\eta_3\mathcal{R}\oplus\eta_4\mathcal{R}$.
Also notice that
$$(1-2u^3)^n=1-2u^3, {~\rm if~} n {~\rm ~is~ odd~ and~} (1-2u^3)^n=1, {~\rm if~} n {~\rm  is~ even.}$$
\noindent For a linear code $\mathcal{C }$ of length $n$ over the ring $\mathcal{R}$, let \vspace{2mm}\\
$\mathcal{C }_1= \{x\in \mathbb{F}_{p}^n : \exists ~y, z, t \in \mathbb{F}_{p}^n {\rm ~such ~that~} \eta_1x+\eta_2y+\eta_3z+\eta_4t \in \mathcal{C }\}$,\vspace{2mm}\\
$\mathcal{C }_2= \{y\in \mathbb{F}_{p}^n : \exists ~x, z, t \in \mathbb{F}_{p}^n {\rm ~such ~that~} \eta_1x+\eta_2y+\eta_3z+\eta_4t \in \mathcal{C }\}$,\vspace{2mm}\\
$\mathcal{C }_3= \{z\in \mathbb{F}_{p}^n : \exists ~x, y, t \in \mathbb{F}_{p}^n {\rm ~such ~that~} \eta_1x+\eta_2y+\eta_3z+\eta_4t \in \mathcal{C }\}$,\vspace{2mm}\\
$\mathcal{C }_4= \{t\in \mathbb{F}_{p}^n : \exists ~x, y, z \in \mathbb{F}_{p}^n {\rm ~such ~that~} \eta_1x+\eta_2y+\eta_3z+\eta_4t \in \mathcal{C }\}$.\vspace{2mm}\\
Then $\mathcal{C}_i, ~i=1,2,3,4$ are linear codes of length $n$ over $\mathbb{F}_{p}$, \vspace{2mm}\\ $\mathcal{C}=\eta_1\mathcal{C}_1\oplus\eta_2\mathcal{C}_2\oplus\eta_3\mathcal{C}_3\oplus\eta_4\mathcal{C}_4
$ and $|\mathcal{C }|= |\mathcal{C }_1|~|\mathcal{C }_2|~|\mathcal{C }_3|~|\mathcal{C }_4|$.\vspace{2mm}

\noindent Let $ \sigma, \gamma, \varrho$ be maps from $\mathcal{R}^n$ to $\mathcal{R}^n$ given by

$$\begin{array}{ll}  \sigma(r_0,r_1,\cdots,r_{n-1})= (r_{n-1},r_0,r_1,\cdots, r_{n-2}),&\vspace{2mm}\\
\gamma(r_0,r_1,\cdots,r_{n-1})= (-r_{n-1},r_0,r_1,\cdots, r_{n-2}),& {\rm and }\vspace{2mm}\\
\varrho(r_0,r_1,\cdots,r_{n-1})= (\lambda r_{n-1},r_0,r_1,\cdots, r_{n-2}),\end{array} $$
\noindent respectively where $\lambda$ is a unit in  $\mathcal{R}$. Let $\mathcal{C}$ be a linear code of length $n$ over $\mathcal{R}$. Then $\mathcal{C}$ is said to be cyclic if $\sigma(\mathcal{C})=\mathcal{C}$, negacyclic if $\gamma(\mathcal{C})=\mathcal{C}$ and $\lambda$-constacyclic if $\varrho(\mathcal{C})=\mathcal{C}$. We will take $\lambda= (1-2u^3)$.\vspace{2mm}

For a code $\mathcal{C}$  over $\mathcal{R}$, the dual code  $\mathcal{C}^\bot$ is defined as $$ \mathcal{C}^\bot =\{x\in \mathcal{R}^n |~ x.y=0 ~ {\rm for ~ all~} y \in \mathcal{C}\}$$ where $x.y$ denotes the Euclidean inner product. $\mathcal{C}$ is called self-dual if $\mathcal{C}=\mathcal{C}^\bot$ and self-orthogonal if $\mathcal{C}\subseteq \mathcal{C}^\bot$.

Two codes $\mathcal{C}_1$ and $\mathcal{C}_2$ of length $n$ over the ring $\mathcal{R}$ are called equivalent or monomially equivalent if there exists a monomial matrix $M$ whose  non-zero scalar entries are some units of the ring $\mathcal{R}$ such that $\mathcal{C}_2=\mathcal{C}_1 M$. More precisely we say that  codes $\mathcal{C}_1$ and $\mathcal{C}_2$ of length $n$ over $\mathcal{R}$ are  equivalent if there are $n$ units $\pi_0,\pi_1,\cdots,\pi_{n-1}$ of $\mathcal{R}$ and a permutation $\tau$ of $n$ coordinate positions such that if $(c_0,c_1,\cdots,c_{n-1}) \in  \mathcal{C}_1$ then $\tau((\pi_0c_0,\pi_1c_1,\cdots,\pi_{n-1}c_{n-1})) \in  \mathcal{C}_2$. If $\pi_i=1$ for each $i$, we say that   $\mathcal{C}_1$ and $\mathcal{C}_2$  are permutation  equivalent. If $\tau$ is identity permutation we say that $\mathcal{C}_1$ and $\mathcal{C}_2$  are scalar  equivalent. Equivalent codes have the same minimum distance.

For a unit $\lambda \in \mathcal{R}$, it is well known that a $\lambda$-constacyclic code of length $n$ over $\mathcal{R}$ can be identified as an ideal in the quotient ring $\frac{\mathcal{R}[x]}{\langle x^{n}-\lambda\rangle}$ via the $\mathcal{R}$-module  isomorphism  $\pi : \mathcal{R}^{n} \rightarrow \frac{\mathcal{R}[x]}{\langle x^{n}-\lambda\rangle}$ given by $$( a_0,a_1,\cdots, a_{n-1}) \mapsto a_0 + a_1x + \cdots + a_{n-1}x^{n-1} ~~({\rm mod~}\langle x^{n}-\lambda \rangle).$$

Define a Gray map $\Phi : \mathcal{R}\rightarrow \mathbb{F}_p^2$ by $a+bu+cu^2+du^3\mapsto(-d,2a+d)$. This map can be extended to $\mathcal{R}^n$ as ~ $\Phi : \mathcal{R}^n\rightarrow\mathbb{F}_p^{2n}$ by $$(r_0,r_1,\cdots,r_{n-1})\mapsto(-d_0,-d_1,\cdots,-d_{n-1},2a_0+d_0,2a_1+d_1,\cdots,2a_{n-1}+d_{n-1}),$$ where $r_i=a_i+b_iu+c_iu^2+d_iu^3, 0\leq i\leq{n-1}$.\vspace{2mm}

\noindent The map $\Phi$ is a $\mathbb{F}_p$-linear map but not one to one.\vspace{2mm}

Define the Gray weight of a an element $r \in \mathcal{R}$ by $w_{G_{\Phi}}(r) =w_H(\Phi(r))$, the Hamming weight of $\Phi(r)$. The Gray weight of  a codeword
$c=(c_0,c_1,\cdots,c_{n-1})\in \mathcal{R}^n$ is defined as $w_{G_{\Phi}}(c)=\sum_{i=0}^{n-1}w_{G_{\Phi}}(c_i)=w_H(\Phi(c))$. For any two elements $c_1, c_2 \in \mathcal{R}^n$, the Gray distance $d_{G_{\Phi}}$ is given by $d_{G_{\Phi}}(c_1,c_2)=w_{G_{\Phi}}(c_1-c_2)=w_H(\Phi(c_1)-\Phi(c_2))$. \vspace{2mm}

\noindent \textbf{Proposition 1.} The Gray map $\Phi$ is a distance preserving map from ($\mathcal{R}^n$, Gray distance) to ($\mathbb{F}_p^{2n}$, Hamming distance).\vspace{2mm}

 Here the Gray image of a self-dual code over $\mathcal{R}$  may not be self-dual as has been for the codes over the ring
$\mathbb{F}_{p}+u\mathbb{F}_{p}, ~u^2=u$, see Proposition 2.5 of [11]. But if $\mathcal{C}$ is a sub code of $(\mathbb{F}_p+u^3\mathbb{F}_p)^n$, we have \vspace{2mm}

\noindent \textbf{Proposition 2.} Let $\mathcal{C}$ be a code of length $n$ over $\mathcal{R}$ such that $\mathcal{C}\subset(\mathbb{F}_p+u^3\mathbb{F}_p)^n$. If $\mathcal{C}$ is self-orthogonal, then so is $\Phi(\mathcal{C}).$\vspace{2mm}\\
 \textbf{Proof.} Let $r_1, r_2\in C$, so that $r_1=a_1+u^3d_1$ and $r_2=a_2+u^3d_2$, where $a_1, a_2, d_1, d_2\in\mathbb{F}_p^n.$ Assume that $\mathcal{C}$ is self-orthogonal. Then $r_1.r_2=0$ gives   $a_1.a_2=0~{\mbox {and}}~a_1.d_2+a_2.d_1+d_1.d_2=0$. So we have
 $$\begin{array}{ll}\Phi(r_1).\Phi(r_2)&=(-d_1,2a_1+d_1).(-d_2,2a_2+d_2)\\&=4a_1.a_2+2a_1.d_2+2a_2.d_1+2d_1.d_2=0,\end{array}$$ which implies that $\Phi(C)$ is self-orthogonal.  ~~~~~~~~~~~~~~~~~~~~~~~~$\Box$\vspace{2mm}

\noindent The polynomial correspondence of the Gray map can be defined as:\\ $~~~~~~~~~~~~~~~~~~~ \Phi:  \mathcal{R}[x]/\langle x^n-\lambda \rangle \rightarrow \mathbb{F}_p[x]/\langle x^{2n}-1 \rangle$    given by $$ \Phi(a(x)+b(x)u+c(x)u^2+d(x)u^3) = ~-d(x)+x^n(2a(x)+d(x)).$$
\section{$(1-2u^3)$-constacyclic codes over $\mathcal{R}$}
Theorems 1-3 and 5, in this section, are generalizations of results of [15]. Proofs are similar so we skip some of them. \vspace{2mm}

\noindent \textbf{Theorem 1.} Let $\varrho$ denote the $(1-2u^3)$- constacyclic shift of $\mathcal{R}^n$ and $\sigma$ the cyclic shift of $\mathbb{F}_p^{2n}$. If $\Phi$ is the Gray Map of $\mathcal{R}^n$ into $\mathbb{F}_p^{2n}$, then $\Phi\varrho=\sigma\Phi$.\vspace{2mm}\\
\noindent \textbf{Corollary 1.} The Gray image $\Phi(\mathcal{C})$ of a $(1-2u^3)$- constacyclic code $\mathcal{C}$ over $\mathcal{R}$ of length $n$ is a cyclic code over  $\mathbb{F}_p$ of length $2n$.\vspace{2mm}\\
This is so because if $\varrho(\mathcal{C})=\mathcal{C}$, then $\Phi(\mathcal{C})= \Phi(\varrho(\mathcal{C}))=\Phi\varrho(\mathcal{C})=\sigma \Phi(\mathcal{C})=\sigma (\Phi(\mathcal{C}))$.\vspace{2mm}

\noindent \textbf{Theorem 2.} Let $\mathcal{C}=\eta_1\mathcal{C}_1\oplus\eta_2\mathcal{C}_2\oplus\eta_3\mathcal{C}_3\oplus\eta_4\mathcal{C}_4$ be a code of length $n$ over $\mathcal{R}$. Then
  $\mathcal{C}$ is a $(1-2u^3)$- constacyclic code of length $n$ over $\mathcal{R}$ if and only if $\mathcal{C}_1$ is cyclic and $\mathcal{C}_2, \mathcal{C}_3, \mathcal{C}_4$ are negacyclic  codes of length $n$ over   $\mathbb{F}_p$. Further  $\mathcal{C}$ is a $(1-2u^3)$- constacyclic code if and only if $\mathcal{C}^\bot$ is a $(1-2u^3)$- constacyclic code.\vspace{2mm}

 \noindent{\bf Proof.} Note that $(1-2u^3)\eta_1=\eta_1$, $(1-2u^3)\eta_i=-\eta_i$, for $2\leq i\leq 4$. \\
For $r=(r_0,r_1,\cdots,r_{n-1})\in C$,  $r_i=\eta_1a_i+\eta_2b_i+\eta_3c_i+\eta_4d_i$, $a=(a_0,a_1,\cdots,a_{n-1})$, $b=(b_0,b_1,\cdots,b_{n-1})$, $c=(c_0,c_1,\cdots,c_{n-1})$ and $d=(d_0,d_1,\cdots,d_{n-1})$, we have $a\in C_1$, $b\in C_2$, $c\in C_3$ and $d\in C_4$. One finds that  $\varrho(r)= \eta_1\sigma(a)+\eta_2\gamma(b)+\eta_3\gamma(c)+\eta_4\gamma(d)$.  Therefore $ \varrho(r)\in\mathcal{C}$
 if and only if
   $\sigma(a)\in \mathcal{C}_1$, $\gamma(b)\in \mathcal{C}_2$, $\gamma(c)\in \mathcal{C}_3$, $\gamma(d)\in \mathcal{C}_4$.
  The second statement holds because  $(1-2u^3)^{-1}=(1-2u^3)$. ~~~~~~~~~~~~$\Box$\vspace{2mm}

 \noindent \textbf{Theorem 3.} Let $n$ be an odd natural number. Then a $(1-2u^3)$- constacyclic code of length $n$ over $\mathcal{R}$ is  equivalent to a cyclic code of length $n$ over $\mathcal{R}$ by the ring isomorphism $\mu : \mathcal{R}[x]/ \langle (x^n-1)\rangle\rightarrow  \mathcal{R}[x]/\langle x^n-(1-2u^3)\rangle$ given by $ \mu(a(x))= a((1-2u^3) x)$. \vspace{2mm}

\noindent This is so because one finds that for $n$ odd, $a(x)\equiv b(x)( {\rm ~mod~}x^n-1 )$ if and only if $a((1-2u^3)x)\equiv b((1-2u^3)x) ({\rm ~mod~}x^n-(1-2u^3))$.\vspace{2mm}

 If $n$ is even, the following Lemma gives conditions on the prime $p$ under which a $(1-2u^3)$- constacyclic code of length $n$ over $\mathcal{R}$ is equivalent to a cyclic code of length $n$ over $\mathcal{R}$.  By Theorem 2.1, Corollary 2.1 and Lemma 2.1 of [1], we have \vspace{2mm}

 \noindent \textbf{Lemma 1.} Let $p$ is an odd prime and $\alpha$ be a primitive element of $\mathbb{F}_p$. Let $\lambda = \alpha^i$, for some $i,~0\le i \le p-2$, Then the following holds: \vspace{2mm}\\
$ \begin{array}{ll}(i)& x^n=\lambda {\rm ~ has~ a ~solution ~in~}  \mathbb{F}_p {\rm ~if~ and~ only~ if~ gcd}(n, p-1)| ~i;\vspace{2mm} \\ (ii)& {\rm If~} \mathbb{F}_p^* {\rm ~contains ~ an~ }n{\rm th ~ root~} \delta {\rm ~ of~} \lambda  {\rm ~ i.e.~}  \delta^n=\lambda, {\rm ~ then~  a~}  \lambda{\rm - constacyclic } \\& {\rm code~ of~ length~} n {\rm ~ over~}  \mathbb{F}_p {\rm ~ is~ scalar~ equivalent~ to~ a ~cyclic~ code~ of~ length~} n  \\& {\rm  over~} \mathbb{F}_p {\rm ~  by~ the~ map~} a(x)\mapsto a(\delta^{-1} x).\end{array}$ \vspace{3mm}

  \noindent \textbf{Theorem 4.} A $(1-2u^3)$- constacyclic code of length $n$ over $\mathcal{R}$ is scalar equivalent to a cyclic code of length $n$ over $\mathcal{R}$ if and only if either $n$ is odd or $n=2^bn',~n'$ odd, $b\geq 1$ and $p\equiv 1 ({\rm ~mod} ~ 2^{b+1})$.\vspace{2mm}

   \noindent  { \bf Proof.} Let $\mathcal{C}=\eta_1\mathcal{C}_1\oplus\eta_2\mathcal{C}_2\oplus\eta_3\mathcal{C}_3\oplus\eta_4\mathcal{C}_4$ be  a
   $(1-2u^3)$- constacyclic code of length $n$ over $\mathcal{R}$ where $\mathcal{C}_1$ is cyclic and $\mathcal{C}_2, \mathcal{C}_3, \mathcal{C}_4$ are negacyclic  codes of length $n$ over   $\mathbb{F}_p$. Now $\mathcal{C}$ is equivalent to a cyclic code if and only if the negacyclic codes $\mathcal{C}_2, \mathcal{C}_3, \mathcal{C}_4$ are equivalent to some cyclic codes. As $-1=\alpha^{\frac{p-1}{2}}$,  an $n$th root of $-1$ exists in $\mathbb{F}_p$ if and only if  ${\rm gcd}(n, p-1)| \frac{p-1}{2}$, by Lemma 1 $(i)$.\\
    If $n$ is odd, we have ${\rm gcd}(n, p-1)= {\rm gcd}(n,\frac{p-1}{2})$ which clearly divides $\frac{p-1}{2}$, $p$ being an odd prime. \\  Let now $n=2^b n',~n'$ odd, $b\geq 1$. If $p\equiv 1 ({\rm ~mod} ~ 2^{b+1})$, write $p-1=2^{b+1}q$. Then ${\rm gcd}(n,p-1)={\rm gcd}(2^bn',2^{b+1}q)=2^b~ {\rm gcd}(n',2q)= 2^b ~{\rm gcd}(n',q)$, $n'$ being odd. So ${\rm gcd}(n,p-1)$  divides $2^bq= \frac{p-1}{2}$. Conversely suppose $p-1=2^r y,~y$ odd, and let ${\rm gcd}(n,p-1)={\rm gcd}(2^bn',2^{r}y)$ divide $\frac{p-1}{2}=2^{r-1}y$. Then we must have $b\leq r-1$ i.e. $r\geq b+1$ which implies $p\equiv 1 ({\rm ~mod} ~ 2^{b+1})$. Now the result follows from Lemma 1.  \vspace{2mm}

    \noindent \textbf{Theorem 5.} Let $\mathcal{C}$ be a $(1-2u^3)$- constacyclic code of length $n$ over $\mathcal{R}$. Then \vspace{2mm}\\ $\begin{array}{ll}
     (i)& \mathcal{C}=\langle\eta_1g_1(x),\eta_2g_2(x),\eta_3g_3(x),\eta_4g_4(x)\rangle {\rm ~where~} g_1(x), g_2(x), g_3(x), g_4(x)\\&{\rm~ are ~the~} {\rm monic ~generator ~polynomials~ of~} \mathcal{C}_1, \mathcal{C}_2, \mathcal{C}_3, \mathcal{C}_4 {\rm ~respectively}. {\rm ~Further~}\\&  g(x)=\eta_1g_1(x)+\eta_2g_2(x)+\eta_3g_3(x)+\eta_4g_4(x){\rm ~is  ~the ~unique ~polynomial~}  \\& {\rm~such ~that~} \mathcal{C}=\langle g(x)\rangle. \vspace{2mm}\\
   (ii)& |\mathcal{C}|=p^{4n-\sum_{i=1}^4deg(g_i)}.\vspace{2mm}\end{array}$\\$\begin{array}{ll}
  (iii)& {\rm Suppose~ that~} g_1(x)h_1(x)=x^n-1 {\rm~ and~} g_i(x)h_i(x)=x^n+1,~ 2\leq i\leq 4.\\&
   {\rm Let ~} h(x)=\eta_1h_1(x)+\eta_2h_2(x)+\eta_3h_3(x)+\eta_4h_4(x), {\rm~then~}\\&
     g(x)h(x)=x^n-(1-2u^3).\vspace{2mm}\end{array}$\\$\begin{array}{ll}(iv)& {\rm If ~gcd} (f_i(x),h_i(x))=1 {\rm~ for~} 1\leq i\leq 4, {\rm ~ then~gcd} (f(x),h(x))=1\\&
     {\rm~ and ~} \mathcal{C}=\langle g(x)f(x)\rangle, {\rm ~where~} f(x)=\eta_1f_1(x)+\eta_2f_2(x)+\eta_3f_3(x)+\eta_4f_4(x) .\vspace{2mm}\end{array}$\\$\begin{array}{ll}(v)& \mathcal{C}^\perp=\eta_1\mathcal{C}_1^\perp\oplus\eta_2\mathcal{C}_2^\perp\oplus\eta_3\mathcal{C}_3^\perp
     \oplus\eta_4\mathcal{C}_4^\perp.\vspace{2mm}\\ (vi) & \mathcal{C}^\perp=\langle h^\perp(x)\rangle, {\rm~~
    where~} h^\perp(x)=\eta_1h_1^\perp(x)+\eta_2h_2^\perp(x)+\eta_3h_3^\perp(x)+\eta_4h_4^\perp(x),\\&
     h_i^\perp(x){\rm ~ is~ the~ reciprocal} {\rm
   ~ polynomial ~of~} h_i(x), ~1\leq i\leq 4 .  \vspace{2mm}\\ (vii)& |\mathcal{C}^\perp|=p^{\sum_{i=1}^4 deg(g_i)}.\vspace{2mm}\\
      (viii) &{\rm Let ~} r_1 {\rm ~ denote ~ the ~ number ~of ~irreducible~factors~of~} x^n-1 {\rm~over~} \mathbb{F}_p {\rm ~and ~} r_2 \\&{\rm  denote ~ the ~ number ~of ~irreducible~factors~of~} x^n+1 {\rm~over~} \mathbb{F}_p {\rm ~then ~the}\\&{\rm number~ of  ~}(1-2u^3){\rm-
      constacyclic ~codes ~of ~length ~} n {\rm ~over ~} \mathcal{R} {\rm ~ is~} 2^{r_1}8^{r_2}. \vspace{2mm}\\ (ix)& \mathcal{C} {\rm ~is~ self ~ dual ~ if ~ and ~only ~if ~} \mathcal{C}_1, \mathcal{C}_2, \mathcal{C}_3, \mathcal{C}_4 {\rm ~ are~self~dual}.
   \end{array}$\\

  \noindent  { \bf Proof.} If $c(x) \in \mathcal{C}=\eta_1\mathcal{C}_1\oplus\eta_2\mathcal{C}_2\oplus\eta_3\mathcal{C}_3\oplus\eta_4\mathcal{C}_4
  $, then $c(x)=\sum_{i=1}^4\eta_ig_i(x)f_i(x)$ where $f_1(x)\in \mathbb{F}_p[x]/(x^n-1)$ and $f_i(x)\in \mathbb{F}_p[x]/(x^n+1)\mbox{~for~} i=2,3,4$. Therefore $\mathcal{C}\subseteq \langle\eta_1g_1(x),\eta_2g_2(x),\eta_3g_3(x),\eta_4g_4(x)\rangle$. Conversely let $f(x)= \eta_1g_1(x)r_1(x)+\eta_2g_2(x)r_2(x)+\eta_3g_3(x)r_3(x)+\eta_4g_4(x)r_4(x)\in \langle\eta_1g_1(x),\eta_2g_2(x),\eta_3g_3(x),\eta_4g_4(x)\rangle $, where $r_i(x) \in \mathcal{R}[x],$ for $i=1,2,3,4$. Find $a_i(x),b_i(x),c_i(x),d_i(x) \in  \mathbb{F}_p[x]$ such that $r
  _i(x)=\eta_1a_i(x)+\eta_2b_i(x)+\eta_3c_i(x)+\eta_4d_i(x)$. Then $f(x)= \eta_1g_1(x)a_1(x)+\eta_2g_2(x)b
  _2(x)+\eta_3g_3(x)c_3(x)+\eta_4g_4(x)d_4(x)\in \mathcal{C}$.
    Since $\eta_ig_i(x)=\eta_ig(x)$, we find that $\mathcal{C}=\langle g(x)\rangle.$ This proves $(i)$ and $(ii)$ is clear.
    \vspace{2mm}

  To prove $(iii)$, one finds that $g(x)h(x)=\eta_1(x^n-1)+\eta_2(x^n+1)+\eta_3(x^n+1)+\eta_4(x^n+1)= x^n-(1-2u^3)$, so that $g(x)|(x^n-(1-2u^3))$.   \vspace{2mm}

      As ${\rm gcd}(f_i(x),h_i(x))=1$ for $1\leq i\leq 4$, there exists $a_i(x), b_i(x)\in \mathcal{R}[x]$ such that $a_i(x)f_i(x)+b_i(x)h_i(x)=1$. Take $a(x)=\sum_{i=1}^4\eta_ia_i(x)$ and $b(x)=\sum_{i=1}^4\eta_ib_i(x)$, then  we find that
     $a(x)f(x)+b(x)h(x)=\sum_{i=1}^4\eta_i(a_i(x)f_i(x)+b_i(x)h_i(x))=1$. Hence ${\rm gcd}(f(x),h(x))=1$. Now  $f(x)g(x)$ has the same zeros as that of $g(x)$, so
       we get $\mathcal{C}=\langle g(x)f(x)\rangle$. This proves $(iv)$.\vspace{2mm}\\
       $(v)$-$(vii)$ can be similarly proved; $(viii)$ and $(ix)$ are clear.\vspace{2mm}

\noindent \textbf{Remark 1:} Part $(ix)$ of previous theorem asserts that a $(1-2u^3)$- constacyclic code of length $n$ over $\mathcal{R}$ is self dual if and only if the cyclic code $\mathcal{C}_1$ and the negacyclic codes $\mathcal{C}_2, \mathcal{C}_3, \mathcal{C}_4 $ over $ \mathbb{F}_p$ are self dual. It is well known (see [9]) that a cyclic code of length $n$ over a finite field $ \mathbb{F}_q$ is self-dual if and only if $n$ is even and $q$ is a power of $2$. Hence a $(1-2u^3)$- constacyclic code of length $n$ over $\mathcal{R}$ with  $p$ a prime, $p\equiv 1 ({\rm~ mod}~ 3)$ can not be self-dual. \vspace{2mm}

 \noindent \textbf{Theorem 6.} Let $\mathcal{C}$ be a $(1-2u^3)$- constacyclic code of length $n$ over $\mathcal{R}$ generated by $g(x)=\eta_1g_1(x)+\eta_2g_2(x)+\eta_3g_3(x)+\eta_4g_4(x)$ where $g_i(x)\in  \mathbb{F}_p[x]$ are  the monic generator polynomials  of $\mathcal{C}_1, \mathcal{C}_2, \mathcal{C}_3, \mathcal{C}_4$ respectively.  Then the Gray image $\Phi(\mathcal{C})$ of $\mathcal{C}$ is a cyclic subcode of $\langle g_1(x)\lambda(x)\rangle$ of length $2n$ where $\lambda(x)= {\rm gcd}( g_2(x),g_3(x),g_4(x))$ in $ \mathbb{F}_p[x]$.

 \noindent{\bf Proof.} Let  $g_1(x)h_1(x)=x^n-1, g_i(x)h_i(x)=x^n+1$ for $i=2,3,4$ where $h_j(x)\in  \mathbb{F}_p[x]$. Let $r(x)\in \mathcal{C}$, so that $r(x)=g(x)f(x)$ where $f(x)=\eta_1f_1(x)+\eta_2f_2(x)+\eta_3f_3(x)+\eta_4f_4(x)\in \mathcal{R}[x]/\langle x^n-(1-2u^3)\rangle $. Then \vspace{2mm} \\
 $\begin{array}{ll} r(x)&= (\eta_1g_1+\eta_2g_2+\eta_3g_3+\eta_4g_4)(\eta_1f_1+\eta_2f_2+\eta_3f_3+\eta_4f_4)\\&
 =\eta_1g_1f_1+\eta_2g_2f_2+\eta_3g_3f_3+\eta_4g_4f_4\\&=g_1f_1+3^{-1}(g_2f_2+\xi g_3f_3+\xi^2 g_4f_4)u+3^{-1}
 (g_2f_2+\xi^2 g_3f_3+\xi g_4f_4)u^2\\&~~~~~~~~+3^{-1}
 (-g_1f_1+g_2f_2+g_3f_3+g_4f_4)u^3.\end{array}$\\
 \noindent Therefore \vspace{2mm} \\
$ \begin{array}{ll}\Phi(r(x))&= g_1f_1-g_2f_2-g_3f_3-g_4f_4+x^n(g_1f_1+g_2f_2+g_3f_3+g_4f_4)\\&=g_1f_1(x^n+1)+g_2f_2(x^n-1)+g_3f_3(x^n-1)+g_4f_4(x^n-1)\\&=
 g_1f_1g_2h_2+g_2f_2g_1h_1+g_3f_3g_1h_1+g_4f_4g_1h_1\\&=g_1(g_2(f_1h_2+f_2h_1)+g_3f_3h_1+g_4f_4h_1)\\&\subseteq \langle g_1g_2, g_1g_3,g_1g_4\rangle \\&= \langle g_1(x)\lambda(x)\rangle. \end{array}$  ~~~~~~~~~~~~~~~~~~~~~~~~$\Box$\vspace{2mm}\\
We give below some examples to illustrate the theory.\vspace{2mm}\\
  \noindent{\bf Example 1:} Let $p=7$ and $n=5$. We have\vspace{2mm}\\
$ \begin{array}{l} x^5-1=(x-1)(x^4+x^3+x^2+x+1)\\x^5+1=(x+1)(x^4-x^3+x^2-x+1)\end{array}$\\
 as the factorization of $x^5-1$ and $x^5+1$ into irreducible factors over $\mathbb{F}_7$. Then the number of  $(1-2u^3)$-constacyclic codes of length $5$ over $\mathcal{R}$ is $2^28^2$.  Let $g_1(x)=x^4+x^3+x^2+x+1$, $g_2(x)=x^4-x^3+x^2-x+1$, $g_3(x)=x+1$, $g_4(x)=x+1$,  so that gcd$(g_2,g_3,g_4)=1$. $\mathcal{C}_1=\langle g_1(x)\rangle$ is the repetition code of length $5$ and $\mathcal{C}_i=\langle g_i(x)\rangle$, $i=2,3,4$ are negacyclic codes of length 5 over $\mathbb{F}_7$. Then $g(x)=\eta_1g_1(x)+\eta_2g_2(x)+\eta_3g_3(x)+\eta_4g_4(x)= x^4+x^3+x^2+x+1+5(u+u^2)(x^4-x^3+x^2-2x)+u^3(4x^4-6x^3+4x^2+4x)$ and $\langle g(x)\rangle$ is a  $(1-2u^3)$- constacyclic code of length $5$ over $\mathcal{R}$. $\mathcal{C}$ is equivalent to a cyclic code generated by  $\eta_1g_1(x)+\eta_2g_2(-x)+\eta_3g_3(-x)+\eta_4g_4(-x)= x^4+x^3+x^2+x+1+5(u+u^2)(x^4+x^3+x^2+3x+6)+u^3(4x^4+4x^3+4x^2+x)$. Further $\Phi(g(x))= -4x^4+6x^3-4x^2-4x+x^5(2x^4+2x^3+2x^2+2x+2+4x^4-6x^3+4x^2+4x)= 6x^9-4x^8+6x^7+6x^6+2x^5-4x^4+6x^3-4x^2-4x$ and $\Phi(\mathcal{C})$ is a cyclic subcode of $\langle x^4+x^3+x^2+x+1 \rangle$ of length 10 over $\mathbb{F}_7$ with minimum distance at least 2. Note that $\Phi(g(x))$ is not the generator polynomial of $\Phi(\mathcal{C})$.\vspace{2mm}

 \noindent{\bf Example 2:} Let $p=7$ and $n=8$. We have
$$ \begin{array}{l} x^8-1=(x-1)(x+1)(x^2+1)(x^2-3x+1)(x^2+3x+1)\\x^8+1=(x^2-x-1)(x^2-3x-1)(x^2+x-1)(x^2+3x-1)\end{array}$$
 as the factorization of $x^8-1$ and $x^8+1$ respectively into irreducible factors over $\mathbb{F}_7$. Then the number of   $(1-2u^3)$-constacyclic codes of length $8$ over $\mathcal{R}$ is $2^58^4$. Let $g_1(x)=x^2-3x+1$, $g_2(x)=x^2-x-1$, $g_3(x)=x^2-3x-1$, $g_4(x)=x^2+x-1$,  so that gcd$(g_2,g_3,g_4)=1$. $\mathcal{C}_1=\langle g_1(x)\rangle$ is a cyclic code of length $8$ and $\mathcal{C}_i=\langle g_i(x)\rangle$, $i=2,3,4$ are negacyclic codes of length 8 over $\mathbb{F}_7$. Then $g(x)=\eta_1g_1(x)+\eta_2g_2(x)+\eta_3g_3(x)+\eta_4g_4(x)= (x^2-3x+1)-ux+u^2x+u^3(2x-2)$ and $\mathcal{C}$ is a  $(1-2u^3)$- constacyclic code of length $8$ over $\mathcal{R}$ generated by $g(x)$. Here $\mathcal{C}$ is not equivalent to a cyclic code. Further $\Phi(g(x))= -2x+2+x^8(2x^2-6x+2+2x-2)= 2x^{10}-4x^9-2x+2$ and $\Phi(\mathcal{C})$ is a cyclic subcode of $\langle x^2-3x+1\rangle$ of length 16 over $\mathbb{F}_7$. This can be generalized to $(1-2u^3)$- constacyclic code of length $2^m, ~m\geq 3$ over $\mathcal{R}$ as\vspace{2mm}\\
 $ x^{2^m}+1
 =(x^{2^{m-2}}-x^{2^{m-3}}-1)(x^{2^{m-2}}-3x^{2^{m-3}}-1)(x^{2^{m-2}}+x^{2^{m-3}}-1)
 (x^{2^{m-2}}+3x^{2^{m-3}}-1)$\vspace{2mm}\\and\\ $x^{2^m}-1=(x-1)(x+1)(x^2+1)(x^2-3x+1)(x^2+3x+1)\times$\\$~~~~~~~~~~~~~~~~~~~{\displaystyle\prod_{j=3}^{m-1}}
 (x^{2^{j-2}}-x^{2^{j-3}}-1)(x^{2^{j-2}}-3x^{2^{j-3}}-1)(x^{2^{j-2}}+x^{2^{j-3}}-1)
 (x^{2^{j-2}}+3x^{2^{j-3}}-1)$\vspace{2mm}
  are the factorizations of $x^{2^m}+1$ and $x^{2^m}-1$ into irreducible factors over $\mathbb{F}_{7}$ (see [2]). \vspace{2mm}

 \noindent{\bf Example 3:} Let $p=19$ and $n=3^m$. We have, following [17]\vspace{2mm}\\
 $ x^{3^m}-1=(x-1)(x-7)(x-11){\displaystyle \prod_{j=0}^{m-2}}\{(x^{3^j}-4)(x^{3^j}-6)(x^{3^j}-9)(x^{3^j}-5)(x^{3^j}+2)(x^{3^j}+3)\}$\vspace{2mm}\\and\\
  $x^{3^m}+1=(x+1)(x+7)(x+11){\displaystyle \prod_{j=0}^{m-2}}\{(x^{3^j}+4)(x^{3^j}+6)(x^{3^j}+9)(x^{3^j}+5)(x^{3^j}-2)(x^{3^j}-3)\}$\vspace{2mm}\\
 as the factorization of $x^{3^m}-1$ and $x^{3^m}+1$  into irreducible polynomials over $\mathbb{F}_{19}$. Then there are $2^{6m-3}8^{6m-3}$ $(1-2u^3)$-constacyclic codes of length $3^m$ over $\mathcal{R}$ which are equivalent to cyclic codes.\vspace{2mm}

 \noindent {\bf Remark 2 :} The generator polynomials of  repeated root cyclic and negacyclic codes over $\mathbb{F}_p$ of various lengths such as $ 3p^s, 2^np^s, \ell^tp^s$, where $\ell $ is an odd  prime different from $p$ have been obtained by several authors in [2,6,7,17]. Theorem 5 immediately gives generator polynomials of $(1-2u^3)$-constacyclic codes of these lengths over the finite non-chain ring $\mathcal{R}$.

\section{Quadratic residue codes over $\mathcal{R}$}

In this section, quadratic residue codes over $\mathcal{R}$ are defined in terms of their idempotent generators. Let $n=q$ be an odd prime such that $p$ is a quadratic residue modulo $q$. Let $Q_q$ and $N_q$ be the sets of quadratic residues and non-residues modulo $q$ respectively. Let
$$ r(x)= {\displaystyle \prod_{r\in Q_q}}(x-\alpha^r), ~~ n(x)= {\displaystyle \prod_{n\in N_q}}(x-\alpha^n)$$
where $\alpha $ is a primitive $q$th root of unity in some extension field of $\mathbb{F}_p$. Following classical notation of [14], let $\mathbb{Q}$, $\mathbb{N}$ be the QR codes generated by $r(x)$ and $n(x)$ and $\tilde{\mathbb{Q}}$, $\tilde{\mathbb{N}}$ be the expurgated QR codes generated by $(x-1)r(x)$ and $(x-1)n(x)$ respectively. we use the notation \vspace{2mm}\\
$~~~~~~j_1(x)=\sum_{i\in Q_q}x^i , ~ j_2(x)=\sum_{i\in N_q}x^i,$\vspace{2mm}\\ $~~~~~~h(x)=1+j_1(x)+j_2(x)=1+x+x^2+\cdots+x^{q-1}=r(x)n(x). $

\noindent If $p>2$,  $q\equiv \pm 1({\rm mod}~4)$ and $p$ is a quadratic residue modulo $q$, then idempotent generators of  $\mathbb{Q}$, $\mathbb{N}$, $\tilde{\mathbb{Q}}$, $\tilde{\mathbb{N}}$ over $\mathbb{F}_p$ are given by, (see [14])\vspace{2mm}\\
$E_q(x)= \frac{1}{2}(1+\frac{1}{q})+\frac{1}{2}(\frac{1}{q}-\frac{1}{\theta})j_1+\frac{1}{2}(\frac{1}{q}+\frac{1}{\theta})j_2,$\vspace{2mm}\\
$E_n(x)= \frac{1}{2}(1+\frac{1}{q})+\frac{1}{2}(\frac{1}{q}-\frac{1}{\theta})j_2+\frac{1}{2}(\frac{1}{q}+\frac{1}{\theta})j_1,$\vspace{2mm}\\
$F_q(x)= \frac{1}{2}(1-\frac{1}{q})-\frac{1}{2}(\frac{1}{q}+\frac{1}{\theta})j_1-\frac{1}{2}(\frac{1}{q}-\frac{1}{\theta})j_2,$\vspace{2mm}\\
$F_n(x)= \frac{1}{2}(1-\frac{1}{q})-\frac{1}{2}(\frac{1}{q}+\frac{1}{\theta})j_2-\frac{1}{2}(\frac{1}{q}-\frac{1}{\theta})j_1,$\vspace{2mm}\\
respectively where $\theta$ denotes Gaussian sum, $\chi(i)$ denotes Legendre symbol, that is
$$ \theta= {\displaystyle \sum_{i=1}^{q-1}}\chi(i)\alpha^i, ~~~~~~\chi(i)= \left\{ \begin{array}{ll}0,& p|i\\1,&i\in Q_q\\-1, & i\in N_q. \end{array}\right.$$
It is known that $\theta^2=-q$ if $q\equiv 3({\rm mod}~4)$ and $\theta^2=q$ if $q\equiv 1({\rm mod}~4)$.\\
For convenience we write $e_1=E_q(x),e_2= E_n(x), \tilde{e_1}=F_q(x),\tilde{e_2}= F_n(x)$.\vspace{2mm}\\
\noindent{\bf Lemma 2}: The four idempotent generators of quadratic residue codes satisfy $e_1+e_2= 1+\frac{1}{q}h$, $\tilde{e_1}+\tilde{e_2}= 1-\frac{1}{q}h$, $e_1-\tilde{e_1}= \frac{1}{q}h$, $e_2-\tilde{e_2}= \frac{1}{q}h$. Further $e_1e_2= \frac{1}{q}h$ and $\tilde{e_1}\tilde{e_2}=0$.\vspace{2mm}\\
\noindent {\bf Proof} A direct simple calculation gives the first four expressions. To see that $e_1e_2= \frac{1}{q}h$, we note that $e_1e_2= e_1(e_1+e_2-1)=e_1(\frac{1}{q}h)=\frac{1}{q}h$ as $e_1$ is multiplicative unity of the QR code $\mathbb{Q}$ and $\frac{1}{q}h$ is divisible by $r(x)$ so it belongs to $\mathbb{Q}$. Again $\tilde{e_1}\tilde{e_2}= \tilde{e_1}(\tilde{e_1}+\tilde{e_2}-1)=\tilde{e_1}(\frac{-1}{q}h)=0$ as $\tilde{e_1} \in
\tilde{\mathbb{Q}}=\langle (x-1)r(x)\rangle$ and $h(x)=r(x)n(x)$, so $\tilde{e}_1h$ is a multiple of $x^q-1$ and hence zero in $\mathcal{R}_q$. \vspace{2mm}\\
\noindent{\bf Lemma 3}: Let $p$ be a prime, $p\equiv 1$(mod 3) and $\eta_i$ be as defined in (1). Then $\eta_1e_i+\eta_2e_j+\eta_3e_k+\eta_4 e_{\ell}$, $\eta_1\tilde{e}_i+\eta_2\tilde{e}_j+\eta_3\tilde{e}_k+\eta_4\tilde{e}_{\ell}$ , are idempotents in the ring $\mathcal{R}_q= \frac{\mathcal{R}[x]}{\langle x^{q}-1\rangle}$ where $e_i, e_j, e_k, e_{\ell}$ are not all equal and $\tilde{e}_i, \tilde{e}_j, \tilde{e}_k, \tilde{e}_{\ell}$ are not all equal for $i,j,k,\ell \in \{1,2\}$.\vspace{2mm}

Working as in Theorem 5, we have \vspace{2mm}

\noindent{\bf Theorem 7}: Let $\mathcal{C}=\eta_1\mathcal{C}_1\oplus\eta_2 \mathcal{C}_2\oplus \eta_3 \mathcal{C}_3\oplus \eta_4\mathcal{C}_4$ be a linear code of length $n$ over $\mathcal{R}$. Then \vspace{2mm}

\noindent (i)~~ $\mathcal{C}$  is cyclic over $\mathcal{R}$ if and only if $\mathcal{C}_i, ~i=1,2,3,4$ are cyclic over $\mathbb{F}_p$.\vspace{2mm}

\noindent (ii) ~~If $\mathcal{C}_i=\langle  g_i (x)\rangle, ~g_i(x)\in \frac{\mathbb{F}_p[x]}{\langle x^{n}-1\rangle}$, $g_i(x)|(x^n-1)$,  then $\mathcal{C}=\langle \eta_1g_1(x),\eta_2g_2(x),\eta_3g_3(x),\eta_4g_4(x)\rangle$\\ $~~~~~~~~ =\langle g(x)\rangle $ where $g(x)= \eta_1g_1+\eta_2g_2+\eta_3g_3+\eta_4g_4$ and $g(x)|(x^{n}-1)$.\vspace{2mm}

\noindent (iii)~~ Further $|\mathcal{C }|=p^{4n-\sum_{i=1}^{4}deg(g_i)}$.\vspace{2mm}

\noindent (v)~~ Suppose that $g_i(x)h_i(x)=x^n-1,~ 1\leq i\leq 4.$ Let $ h(x)=\eta_1h_1(x)+\eta_2h_2(x)+$ \\$~~~~~~~~~\eta_3h_3(x)+\eta_4h_4(x),$ then
     $g(x)h(x)=x^n-1$. \vspace{2mm}

\noindent(vi)~~ $ \mathcal{C}^\perp=\eta_1\mathcal{C}_1^\perp\oplus\eta_2\mathcal{C}_2^\perp\oplus\eta_3\mathcal{C}_3^\perp
     \oplus\eta_4\mathcal{C}_4^\perp.$ \vspace{2mm}

\noindent (vii)~~ $ \mathcal{C}^\perp=\langle h^\perp(x)\rangle,$
    where $ h^\perp(x)=\eta_1h_1^\perp(x)+\eta_2h_2^\perp(x)+\eta_3h_3^\perp(x)+\eta_4h_4^\perp(x)$,
      $ h_i^\perp(x)$ \\ $~~~~~~~~~$is the reciprocal
    polynomial of $h_i(x), ~1\leq i\leq 4.$  \vspace{2mm}

 \noindent(viii)$~~ |\mathcal{C}^\perp|=p^{\sum_{i=1}^4 deg(g_i)}$.\vspace{2mm}

The following is a well known result :\vspace{2mm}

\noindent{\bf Lemma 4}: (i) Let $C$  be a cyclic code of length $n$ over a finite ring $S$ generated by the idempotent $E$ in $S[x]/\langle x^n-1\rangle$ then $C^{\perp}$ is generated by the idempotent $1-E(x^{-1})$.\vspace{2mm}

\noindent (ii) Let $C$ and $D$ be cyclic codes of length $n$ over a finite ring $S$ generated by the idempotents $E_1$ and $E_2$ in $S[x]/\langle x^n-1\rangle$ then $C\cap D$ and $C+D$ are generated by the idempotents $E_1E_2$ and $E_1+E_2-E_1E_2$ respectively.\vspace{2mm}

\noindent Set \vspace{2mm}\\
$Q_1= \langle \eta_1e_2+\eta_2e_1+\eta_3e_1+\eta_4e_1\rangle, ~ Q_2= \langle \eta_1e_1+\eta_2e_2+\eta_3e_1+\eta_4e_1\rangle,$\vspace{2mm}\\
$Q_3= \langle \eta_1e_1+\eta_2e_1+\eta_3e_2+\eta_4e_1\rangle, ~ Q_4= \langle \eta_1e_1+\eta_2e_1+\eta_3e_1+\eta_4e_2\rangle,$\vspace{2mm}\\
$Q_5= \langle \eta_1e_1+\eta_2e_2+\eta_3e_2+\eta_4e_2\rangle, ~ Q_6= \langle \eta_1e_2+\eta_2e_1+\eta_3e_2+\eta_4e_2\rangle,$\vspace{2mm}\\
$Q_7= \langle \eta_1e_2+\eta_2e_2+\eta_3e_1+\eta_4e_2\rangle, ~ Q_8= \langle \eta_1e_2+\eta_2e_2+\eta_3e_2+\eta_4e_1\rangle,$\vspace{2mm}\\
$Q_9=Q_{(1,2)}= \langle \eta_1e_1+\eta_2e_1+\eta_3e_2+\eta_4e_2\rangle, ~ Q_{10}=Q_{(1,3)}= \langle \eta_1e_1+\eta_2e_2+\eta_3e_1+\eta_4e_2\rangle,$\vspace{2mm}\\
$Q_{11}=Q_{(1,4)}= \langle \eta_1e_1+\eta_2e_2+\eta_3e_2+\eta_4e_1\rangle, ~ Q_{12}=Q_{(3,4)}= \langle \eta_1e_2+\eta_2e_2+\eta_3e_1+\eta_4e_1\rangle,$\vspace{2mm}\\
$Q_{13}=Q_{(2,4)}= \langle \eta_1e_2+\eta_2e_1+\eta_3e_2+\eta_4e_1\rangle, ~ Q_{14}=Q_{(2,3)}= \langle \eta_1e_2+\eta_2e_1+\eta_3e_1+\eta_4e_4
\rangle.$\vspace{2mm}\\
$S_1= \langle \eta_1\tilde{e_2}+\eta_2\tilde{e_1}+\eta_3\tilde{e_1}+\eta_4\tilde{e_1}\rangle, ~ S_2= \langle \eta_1\tilde{e_1}+\eta_2\tilde{e_2}+\eta_3\tilde{e_1}+\eta_4\tilde{e_1}\rangle,$\vspace{2mm}\\
$S_3= \langle \eta_1\tilde{e_1}+\eta_2\tilde{e_1}+\eta_3\tilde{e_2}+\eta_4\tilde{e_1}\rangle, ~ S_4= \langle \eta_1\tilde{e_1}+\eta_2\tilde{e_1}+\eta_3\tilde{e_1}+\eta_4\tilde{e_2}\rangle,$\vspace{2mm}\\
$S_5= \langle \eta_1\tilde{e_1}+\eta_2\tilde{e_2}+\eta_3\tilde{e_2}+\eta_4\tilde{e_2}\rangle, ~ S_6= \langle \eta_1\tilde{e_2}+\eta_2\tilde{e_1}+\eta_3\tilde{e_2}+\eta_4\tilde{e_2}\rangle,$\vspace{2mm}\\
$S_7= \langle \eta_1\tilde{e_2}+\eta_2\tilde{e_2}+\eta_3\tilde{e_1}+\eta_4\tilde{e_2}\rangle, ~ S_8= \langle \eta_1\tilde{e_2}+\eta_2\tilde{e_2}+\eta_3\tilde{e_2}+\eta_4\tilde{e_1}\rangle,$\vspace{2mm}\\
$S_9= S_{(1,2)}=\langle \eta_1\tilde{e_1}+\eta_2\tilde{e_1}+\eta_3\tilde{e_2}+\eta_4\tilde{e_2}\rangle, ~ S_{10}= S_{(1,4)}=\langle \eta_1\tilde{e_1}+\eta_2\tilde{e_2}+\eta_3\tilde{e_2}+\eta_4\tilde{e_1}\rangle,$\vspace{2mm}\\
$S_{11}= S_{(1,3)}=\langle \eta_1\tilde{e_1}+\eta_2\tilde{e_2}+\eta_3\tilde{e_1}+\eta_4\tilde{e_2}\rangle, ~ S_{12}=S_{(3,4)}= \langle \eta_1\tilde{e_2}+\eta_2\tilde{e_2}+\eta_3\tilde{e_1}+\eta_4\tilde{e_1}\rangle,$\vspace{2mm}\\
$S_{13}=S_{(2,3)}= \langle \eta_1\tilde{e_2}+\eta_2\tilde{e_1}+\eta_3\tilde{e_1}+\eta_4\tilde{e_2}\rangle, ~ S_{14}= S_{(2,4)}=\langle \eta_1\tilde{e_2}+\eta_2\tilde{e_1}+\eta_3\tilde{e_2}+\eta_4\tilde{e_1}\rangle.$\vspace{2mm}\\
These twenty eight codes are defined as quadratic residue codes of length $q$ over the ring $\mathcal{R}$. \vspace{2mm}\\
\noindent Since $\sum_{i=1}^{4}\eta_i=1,$ we have for $i=1,2,3,4$\\ $Q_i= \langle (1-\eta_i)e_1+\eta_ie_2\rangle$, $Q_{i+4}= \langle \eta_ie_1+(1-\eta_i)e_2\rangle$,\\
$S_i= \langle (1-\eta_i)\tilde{e_1}+\eta_i\tilde{e_2}\rangle$, $S_{i+4}= \langle \eta_i\tilde{e_1}+(1-\eta_i)\tilde{e_2}\rangle$ .\vspace{2mm}

\noindent For $k=9,10,11$ and $(i,j)=(1,2),(1,3),(1,4)$  respectively we have \\
$Q_k= \langle (\eta_i+\eta_j)e_1+(1-\eta_i-\eta_j)e_2\rangle$,
$Q_{k+3}= \langle (\eta_i+\eta_j)e_2+(1-\eta_i-\eta_j)e_1\rangle$,  \\
$S_k= \langle (\eta_i+\eta_j)\tilde{e_1}+(1-\eta_i-\eta_j)\tilde{e_2}\rangle$,
$S_{k+3}= \langle (\eta_i+\eta_j)\tilde{e_2}+(1-\eta_i-\eta_j)\tilde{e_1}\rangle$.\vspace{2mm}

\noindent{\bf Theorem 8 :} If $p\equiv 1({\rm mod}~3)$, $q$ an odd prime and $p$ is a quadratic residue modulo $q$, then
the following assertions hold for quadratic residues codes over $\mathcal{R}$:
\vspace{2mm}\\$\begin{array}{ll}
{\rm (i)} & Q_i {\rm ~is ~equivalent~ to~} Q_{i+4}{\rm ~ and~}  S_i {\rm ~ is~ equivalent~ to~} S_{i+4} {\rm~ for~} i=1,2,3,4, \vspace{2mm}\\
{\rm (ii)}& Q_j {\rm ~is ~equivalent~ to~} Q_{j+3}{\rm ~ and~}  S_j {\rm ~ is~ equivalent~ to~} S_{j+3} {\rm~ for~} j=9,10,11, \vspace{2mm}\\
{\rm (iii)}& Q_i\cap Q_{i+4}= \langle \frac{1}{q}h\rangle, Q_i+ Q_{i+4}= \mathcal{R}_q, {\rm~ for~} i=1,2,3,4,\vspace{2mm}\\& ~~S_i\cap S_{i+4}= \{0\}, S_i+ S_{i+4}= \langle 1-\frac{1}{q}h\rangle {\rm~ for~} i=1,2,3,4,\vspace{2mm}\\
{\rm (iv)}& Q_k\cap Q_{k+3}= \langle \frac{1}{q}h\rangle, Q_k+ Q_{k+3}= \mathcal{R}_q,  k=9,10,11,\vspace{2mm}\\& S_k\cap S_{k+3}= \{0\}, S_k+ S_{k+3}= \langle 1-\frac{1}{q}h\rangle {\rm~ for~} k=9,10,11,\vspace{2mm}\\
{\rm (v)} &S_i\cap \langle \frac{1}{q}h\rangle = \{0\}, S_i+ \langle \frac{1}{q}h\rangle = Q_i {\rm~ for~} i=1,2,\cdots,14,\vspace{2mm}\\
{\rm (vi)}& |Q_i|= p^{2(q+1)}, |S_i|= p^{2(q-1)}{\rm~ for~} i=1,2,\cdots,14. \end{array}$\vspace{2mm}\\

\noindent \textbf{Proof:} Let $ n\in N_q$. Let $\mu_n$ be the multiplier map $\mu_n : \mathbb{F}_p \rightarrow \mathbb{F}_p$ given by $\mu_n(a)=an ({\rm mod~} p)$ and acting on polynomials as  $\mu_n(\sum_i f_ix^i)=\sum_if_ix^{\mu_n(i)}$. Then $\mu_n(j_1)=j_2$ and $\mu_n(j_2)=j_1$. Therefore $\mu_n(e_1)=e_2$, $\mu_n(e_2)=e_1$, $\mu_n(\tilde{e}_1)=\tilde{e}_2$, $\mu_n(\tilde{e}_2)=\tilde{e}_1$ and so $\mu_n((1-\eta_i)e_1+\eta_ie_2) =(1-\eta_i)e_2+\eta_ie_1$, $\mu_n((1-\eta_i)\tilde{e}_1+\eta_i\tilde{e}_2) =(1-\eta_i)\tilde{e}_2+\eta_i\tilde{e}_1$, $\mu_n((\eta_i+\eta_j)e_1+(1-\eta_i-\eta_j)e_2)= (\eta_i+\eta_j)e_2+(1-\eta_i-\eta_j)e_1$ and $\mu_n((\eta_i+\eta_j)\tilde{e}_1+(1-\eta_i-\eta_j)\tilde{e}_2)= (\eta_i+\eta_j)\tilde{e}_2+(1-\eta_i-\eta_j)\tilde{e}_1$. This proves (i) and (ii). \vspace{2mm}

Let $E_i= (1-\eta_i)e_1+\eta_ie_2$ and $E_i'=(1-\eta_i)e_2+\eta_ie_1$, $\tilde{E}_i= (1-\eta_i)\tilde{e}_1+\eta_i\tilde{e}_2$ and ${\tilde{E}_i'}=(1-\eta_i)\tilde{e}_2+\eta_i\tilde{e}_1$. We note that $E_i+ {E'_i}= e_1+e_2$, $E_i{E_i'}= e_1e_2$, $\tilde{E}_i+{{\tilde{E}}_i'}= \tilde{e}_1+\tilde{e}_2$ and $\tilde{E}_i{\tilde{E}_i'}= \tilde{e}_1\tilde{e}_2$. Therefore by Lemmas 2 and 4, $Q_i\cap Q_{i+4}= \langle E_i {E}_i'\rangle= \langle \frac{1}{q}h\rangle$,
and $Q_i+ Q_{i+4}=\langle E_i+E_i'-E_i E_i'\rangle= \langle e_1+e_2-e_1e_2\rangle= \mathcal{R}_q
$;  ~$S_i\cap S_{i+4}=  \langle \tilde{E}_i \tilde {E}_i'\rangle= \langle \tilde{e}_1\tilde{e}_2\rangle =0$,
and $S_i+ S_{i+4}=\langle \tilde{E}_i+\tilde{E}_i'-\tilde{E}_i \tilde{E}_i'\rangle= \langle \tilde{e}_1+\tilde{e}_2-\tilde{e}_1 \tilde{e}_2\rangle= \langle 1-\frac{1}{q}h \rangle$. This proves (iii). The proof of (iv) is similar.\vspace{2mm}

\noindent Note that, using  Lemma 2, we have for $i=1,2,\cdots,8$\\ $\tilde{E}_i (\frac{1}{q}h)=((1-\eta_i)\tilde{e}_1+\eta_i\tilde{e}_2)(\frac{1}{q}h)= ((1-\eta_i)\tilde{e}_1+\eta_i\tilde{e}_2)(1-\tilde{e}_1-\tilde{e}_2)=0$.\\
 Also  $\tilde{E}_i +\frac{1}{q}h$=
$(1-\eta_i)\tilde{e}_1+\eta_i\tilde{e}_2+(1-\eta_i+\eta_i)(\frac{1}{q}h)=(1-\eta_i)(\tilde{e}_1 +\frac{1}{q}h)+\eta_i(\tilde{e}_2 +\frac{1}{q}h)
=(1-\eta_i)e_1+\eta_ie_2$. Therefore $S_i\cap \langle \frac{1}{q}h\rangle = \langle\tilde{E}_i (\frac{1}{q}h)
\rangle= \{0\},$ and $S_i+ \langle \frac{1}{q}h\rangle = \langle\tilde{E}_i +\frac{1}{q}h-\tilde{E}_i(\frac{1}{q}h)\rangle=
\langle(1-\eta_i)e_1+\eta_ie_2\rangle=Q_i$. This proves (v) for $i=1,2,\cdots,8$. For $i=9,10,\cdots,14$ the proof follows on the same lines.\vspace{2mm}

\noindent Finally for $i=1,2,3,4$ we have $|Q_i\cap Q_{i+4}|= | \frac{1}{q}h|= p^4$, it being a repetition code over  $\mathcal{R}$. Therefore
$$ p^{4q}=|\mathcal{R}_q|=|Q_i+ Q_{i+4}|=\frac{|Q_i| |Q_{i+4}|}{|Q_i\cap Q_{i+4}|}=\frac{|Q_i|^2}{p^4}.$$
This gives $|Q_i|=p^{2(q+1)}$. Similar argument gives $|Q_i|=p^{2(q+1)}$ for $i=9,10,11$. Now for each $i,~1\leq i\leq 14$,  we find  that \vspace{2mm}\\
$~~~~~~~~~~~~~~~~~~~~~~~ p^{2(q+1)}=|Q_i|=|S_i+ \langle\frac{1}{q}h\rangle|=|S_i| |\langle\frac{1}{q}h\rangle|=|S_i| p^4.$\\
since $ |S_i\cap \langle\frac{1}{q}h\rangle|= | \langle 0\rangle|= 1$. This gives $|S_i|=p^{2(q-1)}$.  ~~~~~~~~~~~~~~~~~~~~~~~~$\Box$\vspace{2mm}

\noindent{\bf Theorem 9 :} If $p\equiv 1({\rm mod}~3)$, $q\equiv  3({\rm mod}~4)$, $p$ is a quadratic residue modulo $q$, then
the following assertions hold for quadratic residues codes over $\mathcal{R}$.
\vspace{2mm}\\
$\begin{array}{ll}
{\rm (i)} & Q_i^{\perp} =  S_i {\rm~ for~} i=1,2,\cdots,14. \vspace{2mm}\\
{\rm (ii)}& S_i {\rm ~is~ self~ orthogonal ~ for~} i=1,2,\cdots,14. \end{array} $\vspace{2mm}

\noindent \textbf{Proof:} As $q\equiv  3({\rm mod}~4)$, $-1$ is a quadratic nonresidue. Therefore $j_1(x^{-1})=j_2$ and $j_2(x^{-1})=j_1$
and so $1-e_1(x^{-1})=\tilde{e}_1(x)$ and $1-e_2(x^{-1})=\tilde{e}_2(x)$. For $E_i= (1-\eta_i)e_1+\eta_ie_2$, $1-E_i(x^{-1})= 1-(1-\eta_i)e_1(x^{-1})-\eta_ie_2(x^{-1})= 1-(1-\eta_i)(1-\tilde{e}_1(x))-\eta_i(1-\tilde{e}_2(x))=(1-\eta_i)\tilde{e}_1+\eta_i\tilde{e}_2$. Now result (i) follows from Lemma 4. Using (v) of Theorem 8, we have $S_i\subseteq Q_i= S_i^{\perp}$. Therefore $S_i$ are self orthogonal.  ~~~~~~~~~~~~~~~~~~~~~~~~$\Box$\vspace{2mm}

\noindent Similarly we get\vspace{2mm}

\noindent{\bf Theorem 10 :} If $p\equiv 1({\rm mod}~3)$, $q\equiv  1({\rm mod}~4)$, $p$ is a quadratic residue modulo $q$, then
the following assertions hold for quadratic residues codes over $\mathcal{R}$.
\vspace{2mm}\\$ \begin{array}{ll}
{\rm (i)} & Q_i^{\perp} =  S_{i+4}, Q_{i+4}^{\perp} =  S_i {\rm~ for~} i=1,2,3,4. \vspace{2mm}\\
{\rm (ii)}& Q_j^{\perp} =  S_{j+3}, Q_{j+3}^{\perp} =  S_j {\rm~ for~} j=9,10,11. \end{array}$

\section{Extended Quadratic Residue codes over $\mathcal{R}$ and their Gray images}

 The extended QR-codes over $\mathbb{F}_{p}+u\mathbb{F}_{p}+u^2\mathbb{F}_{p}+u^3\mathbb{F}_{p}$ are formed in the same way as the extended QR-codes over $\mathbb{F}_{p}$ are formed. \vspace{2mm}

\noindent{\bf Theorem 11 :} If  $q\equiv  3({\rm mod}~4)$, then the extended QR-codes $\overline{Q_i}$ for $i=1,...,14$ of length $q+1$ are self dual.\vspace{2mm}

 \noindent \textbf{Proof:} Using quadratic reciprocity law, one easily finds that $-q$ is a quadratic residue modulo $p$  in this case. Find an element $r\in \mathbb{F}_{p}$ such that $r^2\equiv -q({\rm~mod~} p)$. As $Q_i=S_i+\langle \frac{1}{q}h\rangle$ for each $i,~i=1,2,\cdots,14$, by Theorem 8, let $\overline{Q_i}$ be the extended QR-code over $\mathcal{R}$ generated by

  $$~~~~~~~\begin{array}{cccccc}
     \infty & 0 & 1 & 2 & \cdots & q-1
   \end{array}\vspace{-2mm}$$ $$ \overline{G_i}=\left(
                  \begin{array}{cccccc}
                    0 &  &  &  &  ~& ~~ \\
                    0 &  &  & G_i & ~ &~~  \\
                    \vdots &  &  &  & ~ & ~ ~\\
                    r & ~1 & 1 & 1 & \cdots ~& 1 ~~\\
                  \end{array}
                \right)$$
 where $G_i$ is a generator matrix for the QR-code $S_i$. The row above the matrix shows the column labeling  by $\mathbb{F}_q\cup \infty$.  Since the all one vector is in $Q_i$ and  $Q_i^{\perp} =  S_{i}$, the last row of $\overline{G_i}$ is orthogonal to all the previous rows of $\overline{G_i}$. The last row is orthogonal to itself also as $r^2=-q $. Further as $S_i$ is self orthogonal by Theorem 9, we find that the code $\overline{Q_i}$ is self orthogonal. Now the result follows from the fact that $|\overline{Q_i}|=p^4 |S_i|=p^{2(q+1)}= |\overline{Q_i}^{\perp}|$.  ~~~~~~~~~~~~~~~~~~~~~~~~$\Box$\vspace{2mm}

\noindent{\bf Theorem 12 :} If  $q\equiv 1({\rm mod}~4)$, then  $\overline{Q_i}^{\perp}=\overline{Q_{i+4}}$ for $i=1,2,3,4$ and $\overline{Q_i}^{\perp}=\overline{Q_{i+3}}$ for $i=9,10,11$.\vspace{2mm}

 \noindent \textbf{Proof:} We give a proof for $i=1,2,3,4$ only; the other part is similar. As $Q_i=S_i+\langle \frac{1}{q}h\rangle$ for each $i$, let $\overline{Q_i}$ be the extended QR-code over $\mathcal{R}$ generated by
 $$~~~~~~~\begin{array}{cccccc}
     \infty & 0 & 1 & 2 & \cdots & q-1
   \end{array}\vspace{-2mm}$$ $$ \overline{G_i}=\left(
                  \begin{array}{cccccc}
                    0 &  &  &  &  ~& ~~ \\
                    0 &  &  & G_i & ~ &~~  \\
                    \vdots &  &  &  & ~ & ~ ~\\
                    r_i & ~1 & 1 & 1 & \cdots ~& 1 ~~\\
                  \end{array}
                \right)$$
 where $r_i=1$ for $i=1,2,3, 4$, $r_i=-q$ for $i=5,6,7,8$ and $G_i$ is a generator matrix for the  QR-code $S_i$.  Since the all one vector is in $Q_i$ and $Q_{i+4}$ both  and  $Q_i^{\perp} =  S_{i+4}$, $Q_{i+4}^{\perp} =  S_{i}$, by Theorem 10, it is orthogonal to all the rows of $G_i$ and $G_{i+4}$. Also the last row of $\overline{G_i}$ is orthogonal to the last row of  $\overline{G_{i+4}}$. Further rows of $G_{i+4}$ are in $S_{i+4}=Q_i^{\perp}$ implies that rows of $G_{i+4}$ are orthogonal to vectors of $Q_i$ and hence to the rows of $G_i$ as $S_i \subseteq Q_i$. Therefore all the rows of $\overline{G_{i+4}}$ are orthogonal to all the rows of $\overline{G_i}$ implying that $\overline{Q_{i+4}} \subseteq \overline{Q_i}^{\perp}$.  Now the result follows from comparing their orders.  ~~~~~~~~~~~~~~~~~~~~~~~~$\Box$\vspace{2mm}

\noindent We define another Gray map $\Psi : \mathcal{R}\rightarrow \mathbb{F}_p^4$ ~~ given by $$\begin{array}{ll}r=a+bu+cu^2+du^3\mapsto&(a,a+b+c+d,a+b\xi^2+c\xi+d,a+b\xi+c\xi^2+d)\vspace{2mm}\\&= (a,b,c,d)\left(
                                                                             \begin{array}{cccc}
                                                                               1 & 1 & 1 & 1 \\
                                                                               0 & 1 & \xi^2 & \xi \\
                                                                               0 & 1 & \xi & \xi^2 \\
                                                                               0 & 1 & 1 & 1 \\
                                                                             \end{array}
                                                                           \right) =(a,b,c,d)M\end{array}
$$ where $\xi \in \mathbb{F}_p $ is as defined in Section 2 and $M$ is a nonsingular matrix of determinant $ 3\xi(1-\xi)$.   This map can be extended to $\mathcal{R}^n$ component wise.

\noindent Corresponding to this Gray map $\Psi$,  the Gray weight of  an element $r \in \mathcal{R}$ is defined as $w_{G_{\psi}}(r) =w_H(\Psi(r))$, the Hamming weight. The  Gray weight $w_{G_{\psi}}(r)$ of  a codeword
$r=(r_0,r_1,\cdots,r_{n-1})\in \mathcal{R}^n$ and the Gray distance $d_{G_{\psi}}$ of two codewords are defined accordingly.\vspace{2mm}


\noindent \textbf{Theorem 13.} The Gray map $\Psi$ is an  $\mathbb{F}_p$ - linear one to one and onto map.  The Gray image $\Psi(\mathcal{C})$ of a self-dual code $\mathcal{C}$  over $\mathcal{R}$ is a self-dual code in $\mathbb{F}_p^{4n}$. It is also distance preserving map from ($\mathcal{R}^n$, Gray distance $d_{G_{\psi}}$) to ($\mathbb{F}_p^{4n}$, Hamming distance).\vspace{2mm}

\noindent \textbf{Proof.} The first assertion holds as $M$ is an invertible matrix over $\mathbb{F}_p$. Let $\mathcal{C}$ be a self dual code over $\mathcal{R}$. Let $r=(r_0,r_1,\cdots,r_{n-1}), s=(s_0,s_1,\cdots,s_{n-1}) \in \mathcal{C}$ where $r_i=a_i+b_iu+c_iu^2+d_iu^3$ and $s_i=a'_i+b'_iu+c'_iu^2+d'_iu^3$. Then $0=r\cdot s=\sum_{i=0}^{n-1}r_is_i$ gives
$$ \begin{array}{l} {\displaystyle \sum_{i=0}^{n-1}}a_ia'_i=0,\\  {\displaystyle \sum_{i=0}^{n-1}}(a_id'_i+a'_id_i+b_ic'_i+b'_ic_i+d_id'_i)=0. \end{array}$$
A simple calculation shows that $$\begin{array}{l}\Psi(r)\cdot \Psi(s)={\displaystyle \sum_{i=0}^{n-1}}\Psi(r_i)\cdot \Psi(s_i)\\ ~~~~~~~~~~~~~~=4{\displaystyle\sum_{i=0}^{n-1}}a_ia'_i+3{\displaystyle\sum_{i=0}^{n-1}}(a_id'_i+a'_id_i+b_ic'_i+b'_ic_i+d_id'_i) =0,\end{array}$$
as $1+\xi+\xi^2=0$. This proves the result.          ~~~~~~~~~~~~~~~~~~~~~~~~$\Box$\vspace{2mm}

 Note that for  a $(1-2u^3)$- constacyclic code $\mathcal{C}$  over $\mathcal{R}$, the image $\Psi(\mathcal{C})$ may not be cyclic, it is only a $\mathbb{F}_p$-linear code,  whereas the image  $\Phi(\mathcal{C})$ was a cyclic code over $\mathbb{F}_p$.\vspace{2mm}

Since $S_i$  is equivalent to $S_{i+4}$  for $ i=1,2,3,4$,
  they have the same weight enumerator. Further $w_{G_{\psi}}(1)= w_H((1,1,1,1))=4=w_H((-q,-q,-q,-q))=w_{G_{\psi}}(-q)$, therefore $\overline{Q_i}$ and $\overline{Q_{i+4}}$ have the same weight enumerator. Similarly $\overline{Q_i}$ and $\overline{Q_{i+3}}$ have the same weight enumerator for $i=9,10,11$.\vspace{2mm}

\noindent \textbf{Corollary 2.} If  $q\equiv  3({\rm mod}~4)$, the Gray images of extended QR-codes $\overline{Q_i}$ i.e. $\Psi(\overline{Q_i})$ for $i=1,2,\cdots, 14$ are self dual codes of length $4(q+1)$ over $\mathbb{F}_p$. If  $q\equiv  1({\rm mod}~4)$, $\Psi(\overline{Q_i})$ are formally self dual codes of length $4(q+1)$ over $\mathbb{F}_p$.\vspace{2mm}

\noindent {\bf Example 4 :} Take $p=7$ and $q=19$ so that $p$ is a quadratic residue modulo $q$.  The sets of quadratic residues and non residues modulo $19$ are $Q_{19}=\{ 1,4,5,6,7,9,11,16,17\}, N_{19}=\{2,3,8,10,12,13,14,15,18\}$. Therefore $j_1(x)=x+x^4+x^5+x^6+x^7+x^9+x^{11}+x^{16}+x^{17}$, $j_2(x)=x^2+x^3+x^8+x^{10}+x^{12}+x^{13}+x^{14}+x^{15}+x^{18}$. $\theta =\alpha -\alpha^2-\alpha^3+\alpha^4+\alpha^5+\alpha^6+\alpha^7-\alpha^8+\alpha^9-\alpha^{10}+\alpha^{11}-\alpha^{12}-\alpha^{13}-\alpha^{14}
-\alpha^{15}+\alpha^{16}+\alpha^{17}-\alpha^{18} =4$ where $\alpha$ is a $19$th root of unity in some extension field of $\mathbb{F}_7$.
$e_1=2+4j_1(x)+6j_2(x)$ and $e_2=2+6j_1(x)+4j_2(x)$. The idempotent generator of QR-code $Q_1$ over $\mathbb{F}_{7}+u\mathbb{F}_{7}+u^2\mathbb{F}_{7}+u^3\mathbb{F}_{7}$ is given by $(1-\eta_1)e_1+\eta_1e_2= 2+6j_1(x)+4j_2(x)+ 2(j_2-j_1)u^3$. $Q_1^{\perp} =  S_1$, $S_1\subseteq S_1^{\perp}$ where $S_1= \langle -1-6j_2(x)-4j_1(x)+ 2(j_2-j_1)u^3\rangle$. Further $\Psi(\overline{Q_1})$ is a self dual code of length $80
$ over $\mathbb{F}_7$.\vspace{2mm}

\noindent \textbf{Acknowledgements.} The authors are very grateful to the anonymous referees for their  comments and suggestions which helped  significantly in improving  the paper.

\end{document}